\documentclass[12pt]{article}
\usepackage{amsmath,amssymb,amsfonts}
\usepackage{color}

\numberwithin{equation}{section}
\newtheorem{thm}{Theorem}

\newtheorem{cor}[thm]{Corollary}
\newtheorem{defn}{Definition}
\def\cwedge{\bigcirc\kern-1.07em\wedge\ }

\begin{document}

\begin{center}
{\LARGE \bf Remarks on $\eta$-Einstein unit tangent bundles}
\end{center}\begin{center}
{\large \bf Y. D. Chai$^{1}$, S. H. Chun$^{1}$, J. H. Park$^{1}$,
and K. Sekigawa$^{2}$}\end{center}

\begin{center}
$^{1}$Sungkyunkwan University,
 Suwon, Korea\\
$^{2}$Niigata University,
    Niigata, Japan\end{center}

\begin{abstract}
We study the geometric properties of the base manifold for the unit
tangent bundle satisfying the $\eta$-Einstein condition with the
standard contact metric structure. One of the main theorems is that
the unit tangent bundle of 4-dimensional Einstein manifold, equipped
with the canonical contact metric structure, is $\eta$-Einstein
manifold if and only if base manifold is the space of constant
sectional curvature 1 or 2.
\end{abstract}



\begin{center}
\section{Introduction}\label{sec1}
\end{center}
\quad\; We consider the $\eta$-Einstein condition, which is suitable
for contact metric manifold in general, that is, the Ricci tensor is
of the form $\rho(X,Y)=\alpha \,g(X,Y) +\beta\,\eta(X)\eta(Y)$ with
$\alpha$ and $\beta$ being smooth functions. In \cite {BoeckxVan2},
Boeckx and Vanchecke determined the unit tangent bundles which are
Einstein with respect to the canonical contact metric structure. In
the present paper, we shall extend their result to the
$\eta$-Einstein case. The scalar curvature of an $\eta$-Einstein
contact metric manifold is not necessarily constant in general,
however, for some special $\eta$-Einstein contact metric manifolds,
we may expect the scalar curvature to be constant. The main theorems
are the following :
\begin{thm}\label{thm4}
Let M be an $n (\geq 2)$-dimensional Riemannian manifold and $T_1 M$
be the unit tangent bundle of $M$ equipped with the canonical
contact metric structure. If $T_1 M$ is an $\eta$-Einstein manifold,
then $\alpha$ and $\beta$ are both constant valued ones on $T_1M$.
\end{thm}

Let $\tau$ be the scalar curvature of $M$, $\rho$ be the Ricci
curvature tensor of $M$, $R$ be the Riemann curvature tensor of $M$
and $\bar\tau$ be the scalar curvature of $T_1 M$. Then we have the
following theorems.

\begin{thm}\label{thm2}
Let $(T_1 M, \eta, \bar g, \phi, \xi)$ be an $\eta$-Einstein
manifold. Then $\tau$, $|\rho|^2$, $|R|^2$, and $\bar\tau$ are all
constant.
\end{thm}

\begin{thm}\label{thm7}
Let $M$ be a 4-dimensional Einstein manifold and $(T_1 M, \eta, \bar
g, \phi, \xi)$ be the unit tangent bundle of $M$ equipped with the
canonical contact metric structure. If $T_1 M$ is an $\eta$-Einstein
manifold if and only if $(M, g)$ is the space of constant sectional
curvature 1 or 2.
\end{thm}

\noindent {\bf Question 1} Can we extend the above Theorem
\ref{thm7} to
higher dimensional cases?\\

 From our arguments in the present paper,
the following question will naturally arise:\\

\noindent {\bf Question 2}  Does there exist $n(\geq4)$ dimensional
Riemannian manifold which is not a space of constant sectional
curvature $1$ or $n-2$, whose unit tangent bundle is
$\eta$-Einstein?\\

In the last section, we consider $\eta$-Einstein unit tangent
bundles of some special base Riemannian manifolds.

\begin{center}
\section{Unit tangent bundle with contact metric structure}\label{sec2}
\end{center}

\quad\; First, we give some preliminaries on a contact metric
manifold. We refer to \cite{Blair1} for more details. A
differentiable $(2n-1)$-dimensional manifold $\bar M$ is said to be
a {\it contact manifold} if it admits a global $1$-form $\eta$ such
that $\eta\wedge(d\eta)^{n-1}\neq 0$ everywhere on $\bar M$, where
the exponent denotes the $(n-1)$-th exterior power. We call such
$\eta$ a {\it contact form} of $\bar M$. It is well known that given
a contact form $\eta$, there exists a unique vector field $\xi$,
which is called the {\it characteristic vector field}, satisfying
$\eta(\xi)=1$ and $d\eta(\xi,\bar X)=0$ for any vector field $\bar
X$ on $\bar M$. A Riemannian metric $\bar g$ is an associated metric
to a contact form $\eta$ if there exists a $(1,1)$-tensor field
$\phi$ satisfying
\begin{equation}\label{2.1}
        \eta(\bar X)=\bar g(\bar X,\xi),\quad d\eta(\bar X,\bar Y)=\bar g(\bar X,\phi \bar Y),
        \quad \phi^2 \bar X=-\bar X+\eta(\bar X)\xi
\end{equation}
where $\bar X$ and $\bar Y$ are vector fields on $\bar M$. From
(\ref{2.1}) it follows that
\begin{equation*}
        \phi\xi=0,\quad \eta\circ\phi=0,
        \quad \bar g(\phi \bar X,\phi \bar Y)=\bar g(\bar X,\bar Y)-\eta(\bar X)\eta(\bar Y).
\end{equation*}
A Riemannian manifold $\bar M$ equipped with structure tensors
$(\eta,\bar g,\phi,\xi)$ satisfying (\ref{2.1}) is said to be a {\it
contact metric manifold}. We assume that a contact metric manifold
$\bar M =(\bar M, \eta,\bar g,\phi,\xi)$ is always oriented by the
$(2n-1)$-form $\eta\wedge(d\eta)^{n-1}$. We denote by $dV$ the
volume form of $\bar M$ with respect to the metric $\bar g$. Then we
may easily observe that $dV=C\eta\wedge(d\eta)^{n-1}$, where
$C=\frac{1}{(n-1)!}$. We now review some elementary facts in a
contact metric manifold. First, for the characteristic vector field
$\xi$, $L_\xi \eta =0$ follows from $\eta(\xi)=1$, $d\eta(\bar
X,\bar Y)=\bar g(\bar X,\phi \bar Y)$ and $d\eta(\xi,\bar X)=0$.
Here $L$ denotes Lie derivation. Next, since $d\circ L_\xi
=L_\xi\circ d$, by taking account of $L_\xi \eta =0$, we have
\begin{equation}\label{4-9}
\begin{aligned}
L_\xi dV &=C L_\xi (\eta\wedge(d\eta)^{n-1})\\
&=C(L_\xi \eta)\wedge(d\eta)^{n-1}\\
&\quad +C\eta\wedge(L_\xi d\eta)\wedge d\eta \wedge\cdots\wedge d\eta\\
&\quad + \cdots +C\eta\wedge d\eta\wedge\cdots\wedge(L_\xi d\eta)\\
&=C\eta\wedge d(L_\xi\eta)\wedge d\eta \wedge\cdots\wedge d\eta\\
&\quad + \cdots +C\eta\wedge d\eta\wedge\cdots\wedge d(L_\xi\eta)\\
&=0.
\end{aligned}
\end{equation}
Since $L_\xi dV =(div\xi)dV$, by the definition of the divergence
$div \xi$ with respect to $dV$ and by (\ref{4-9}), we have
\begin{equation}\label{4-10}
div \xi=0\quad(i.e., \bar\nabla_i \xi^i =0).
\end{equation}
Since $\bar\nabla_{\bar X} \xi$ is orthogonal to $\xi$, we have
immediately
\begin{equation}\label{4-11}
(\bar\nabla_{\bar X} \eta)\xi=0
\end{equation}
for any vector field $\bar X$ on $\bar M$.\vskip.2cm

Let $(M,g)$ be an $n$-dimensional Riemannian manifold and $\nabla$
the associated Levi Civita connection. Its Riemann curvature tensor
$R$ is defined by $R(X,Y)Z=\nabla_X\nabla_Y Z- \nabla_Y\nabla_X Z
-\nabla_{[X,Y]} Z$ for all vector fields $X,Y$ and $Z$ on $M$. The
tangent bundle of $(M,g)$ is denoted by $TM$ and consists of pairs
$(p,u)$, where $p$ is a point in $M$ and $u$ a tangent vector to $M$
at $P$. The mapping $\pi:TM\rightarrow M,\ \pi(p,u)=p$ is the
natural projection from $TM$ onto $M$.
For a vector field $X$ on $M$, its
{\it vertical lift} $X^v$ on $TM$ is the vector field defined by
$X^v \omega =\omega(X)\circ \pi$, where $\omega$ is a 1-form on $M$.
For a Levi Civita connection $\nabla$ on $M$, the {\it horizontal
lift} $X^h$ of $X$ is defined by $X^h \omega =\nabla_X\omega$. The
tangent bundle $TM$ can be endowed in a natural way with a
Riemannian metric $\tilde g$, the so-called {\it Sasaki metric},
depending only on the Riemannian metric $g$ on $M$. It is determined
by
$$
  \tilde g(X^h,Y^h)=\tilde g(X^v,Y^v)=g(X,Y)\circ\pi,\quad
  \tilde g(X^h,Y^v)=0
$$
for all vector fields $X$ and $Y$ on $M$. Also, $TM$ admits an
almost complex structure tensor $J$ defined by $JX^h=X^v$ and
$JX^v=-X^h$. Then $g$ is Hermitian metric for the almost complex
structure $J$. We note that $J$ is integrable if and only if $(M,g)$
is locally flat (\cite{Dombrowski}).

The unit tangent bundle
$\bar\pi:T_1M\rightarrow M $ is a hypersurface of $TM$ given by $g_p
(u,u)=1$. Note that $\bar\pi =\pi\circ i$, where $i$ is the
immersion. A unit normal vector $N=u^v$ to $T_1M$ is given by the
vertical lift of $u$ for $(p,u)$. The horizontal lift of a vector
is tangent to $T_1M$, but the vertical lift of vector is not tangent
to $T_1M$ in general. So, we define the {\it tangential lift} of $X$
to $(p,u)\in T_1 M$ by
$$
  X^t_{(p,u)}=(X-g(X,u)u)^v.
$$
Clearly, the tangent space $T_{(p,u)}T_1 M$ is spanned by vectors
of the form $X^h$ and $X^t$, where $X\in T_p M$.

\vskip.2cm

We now define the standard contact metric structure of the unit tangent
bundle $T_1M$ of a Riemannian manifold $(M,g)$. The metric $g'$ on $T_1M$ is
induced from the Sasaki metric $\tilde g$ on $TM$.
Using the almost complex structure $J$ on $TM$, we define a unit vector field
$\xi'$, a 1-form $\eta'$ and a (1,1)-tensor field $\phi'$ on $T_1M$ by
\begin{equation*}
  \xi'=-JN,\quad \phi'=J-\eta'\otimes N.
\end{equation*}
Since $g'(X,\phi' Y)=2d\eta'(X,Y)$, $(\eta',g',\phi',\xi')$ is not a
contact metric structure. If we rescale by
$$
\xi=2\xi',\quad \eta=\frac{1}{2}\,\eta',\quad \phi=\phi',\quad \bar g=\frac{1}{4}\, g',
$$
we get the standard contact metric structure $(\eta,\bar g,\phi,\xi)$.
These tensors are given by
\begin{equation}\label{3.1}
 \begin{split}
   &\xi=2u^h,\\
   &\phi X^t=-X^h+\frac{1}{2}\, g(X,u)\xi,\quad
   \phi X^h=X^t,\\
   &\eta(X^t)=0,\quad \eta(X^h)=\frac{1}{2}g(X,u),\\
   &\bar g(X^t,Y^t)=\frac{1}{4}(g(X,Y)-g(X,u)g(Y,u)),\\
   &\bar g(X^t,Y^h)=0,\\
   &\bar g(X^h,Y^h)=\frac{1}{4}g(X,Y),
 \end{split}
\end{equation}
where $X$ and $Y$ are vector fields on $M$. From now on, we
consider $T_1 M=(T_1 M,\eta,\bar g,\phi,\xi)$ with the standard contact
metric structure.

\noindent The Levi Civita connection $\bar\nabla$ of $T_1 M$ is
described by
\begin{equation}\label{3.2}
  \begin{split}
    \bar\nabla_{X^t}Y^t&=-g(Y,u)X^t,\\
    \bar\nabla_{X^t}Y^h&=\frac{1}{2}\,(R(u,X)Y)^h,\\
    \bar\nabla_{X^h}Y^t&=(\nabla_{X}Y)^t+\frac{1}{2}\,(R(u,Y)X)^h,\\
    \bar\nabla_{X^h}Y^h&=(\nabla_{X}Y)^h-\frac{1}{2}\,(R(X,Y)u)^t
  \end{split}
\end{equation}
for all vector fields $X$ and $Y$ on $M$.

\noindent Also the Riemann curvature tensor $\bar R$ of $T_1 M$
is given by
\begin{equation}\label{3.3}
  \begin{split}
    \bar R(X^t,Y^t)Z^t&=-(g(X,Z)-g(X,u)g(Z,u))Y^t +(g(Y,Z)-g(Y,u)g(Z,u))X^t,\\
    \bar R(X^t,Y^t)Z^h&=\big\{R(X-g(X,u)u,Y-g(Y,u)u)Z\big\}^h\\
     &\quad+\frac{1}{4}\big\{[R(u,X),R(u,Y)]Z\big\}^h,\\
    \bar R(X^h,Y^t)Z^t&=-\frac{1}{2}\big\{R(Y-g(Y,u)u,Z-g(Z,u)u)X\}^h\\
     &\quad-\frac{1}{4}\{R(u,Y)R(u,Z)X\big\}^h,\\
    \bar R(X^h,Y^t)Z^h&=\frac{1}{2}\big\{R(X,Z)(Y-g(Y,u)u)\}^t-
     \frac{1}{4}\big\{R(X,R(u,Y)Z)u\big\}^t\\
     &\quad+\frac{1}{2}\big\{(\nabla_X R)(u,Y)Z\big\}^h,\\
    \bar R(X^h,Y^h)Z^t&=\big\{R(X,Y)(Z-g(Z,u)u)\big\}^t\\
     &\quad+\frac{1}{4}\big\{R(Y,R(u,Z)X)u-R(X,R(u,Z)Y)u\big\}^t\\
     &\quad+\frac{1}{2}\big\{(\nabla_X R)(u,Z)Y-(\nabla_Y R)(u,Z)X\big\}^h,\\
    \bar R(X^h,Y^h)Z^h&=(R(X,Y)Z)^h+\frac{1}{2}\big\{R(u,R(X,Y)u)Z\big\}^h\\
     &\quad-\frac{1}{4}\big\{R(u,R(Y,Z)u)X-R(u,R(X,Z)u)Y\big\}^h\\
     &\quad+\frac{1}{2}\big\{(\nabla_Z R)(X,Y)u\big\}^t
  \end{split}
\end{equation}
for all vector fields $X$, $Y$ and $Z$ on $M$.

Next, to calculate the Ricci tensor $\bar\rho$ of $T_1 M$
at the point $(p,u)\in T_1M$, let $e_1,\cdots,e_n
=u$ be an orthonormal basis of $T_p M$. Then
$2e_1^t,\cdots,2e_{n-1}^t,$ $2e_1^h,\cdots,2e_n^h =\xi,$ is an
orthonormal basis for ${T_{(p,u)}}T_1M$ and $\bar\rho$ is given by
\begin{equation}\label{3.4}
  \begin{split}
    \bar\rho(X^t,Y^t)&=(n-2)(g(X,Y)-g(X,u)g(Y,u))+
    \frac{1}{4}\sum_{i=1}^{n}g(R(u,X)e_i , R(u,Y)e_i),\\
    \bar\rho(X^t,Y^h)&=\frac{1}{2}((\nabla_u \rho)(X,Y)-(\nabla_X \rho)(u,Y)),\\
    \bar\rho(X^h,Y^h)&=\rho(X,Y)-\frac{1}{2}\sum_{i=1}^{n}g(R(u,e_i)X , R(u,e_i)Y),\\
  \end{split}
\end{equation}
where $\rho$ denotes the Ricci curvature tensor of $M$.
 From this, the scalar curvature $\bar\tau$ is given by
\begin{equation}\label{3.5}
    \bar\tau
    =\tau + (n-1)(n-2)-\frac{1}{4}\sum_{i,j=1}^{n}g(R(u,e_i)e_j, R(u,e_i)e_j),
\end{equation}
where $\tau$ is the scalar curvature of $M$.

\section{Unit tangent bundles with $\eta$-Einstein structure}\label{sec3}

\quad\; We shall introduce the definition of $\eta$-Einstein
manifold.
\begin{defn}
If the Ricci tensor $\bar\rho$ of a contact metric manifold
$(\bar M, \eta, \bar g, \phi, \xi)$ is of the form
$$
\bar\rho(\bar X,\bar Y)=\alpha\,\bar g(\bar X,\bar Y)
+\beta\,\eta(\bar X)\eta(\bar Y)
$$
for smooth functions $\alpha$ and $\beta$, then $\bar M$ is called an {\it
$\eta$-Einstein manifold}.
\end{defn}

Now, let $M =(M,g)$ be a Riemannian manifold and $(T_1 M, \eta, \bar
g, \phi, \xi )$ be the unit tangent bundle of $(M, g)$ equipped with
the canonical contact metric structure $(\eta, \bar g, \phi, \xi)$
defined as in section \ref{sec2}. Take the $\phi$-basis $\{\bar e_i,
\bar e_{i^*}=\phi \bar e_i, \xi=\bar e_*\}$ on $T_1M$. Then the
Ricci tensor $\bar\rho$ with respect to the $\phi$-basis should be
\begin{equation}
\bar\rho
 =\begin{pmatrix}
  \bar\rho(\bar e_i, \bar e_j) & \bar\rho(\bar e_i, \bar e_{j^*}) & \bar\rho(\bar e_i, \bar e_*)\\
  \bar\rho(\bar e_{i^*}, \bar e_j) & \bar\rho(\bar e_{i^*}, \bar e_{j^*}) & \bar\rho(\bar e_{i^*}, \bar e_*)\\
  \bar\rho(\bar e_*, \bar e_j) & \bar\rho(\bar e_*, \bar e_{j^*}) & \bar\rho(\bar e_*, \bar e_*)
 \end{pmatrix}.
\end{equation}
In particular, if $T_1M$ is $\eta$-Einstein, by the definition, the
Ricci tensor $\bar\rho$ is given by
\begin{equation}
\bar\rho
 =\begin{pmatrix}
  \alpha & o      & \cdots & 0      & 0      \\
  0      & \alpha & \cdots & 0      & 0      \\
  \vdots & \vdots &        & \vdots & \vdots \\
  0      & 0      & \cdots & \alpha & 0      \\
  0      & 0      & \cdots & 0      & \alpha+\beta
 \end{pmatrix}
\end{equation}
for some two smooth functions $\alpha$ and $\beta$ on $T_1M$. From
(\ref{3.4}), we have the following theorem.

\begin{thm}\label{thm1}
Let $M$ be an $n$-dimensional Riemannian manifold. Then $T_1 M$ is
$\eta$-Einstein if and only if
\begin{eqnarray}
\label{4.1} && \sum_{i=1}^{n}g(R(u,X)e_i , R(u,Y)e_i)
                  =(\alpha-4n+8)(g(X,Y)-g(X,u)g(Y,u)),\\
\label{4.2} && (\nabla_u \rho)(X,Y)=(\nabla_X \rho)(u,Y),\\
\label{4.3} && \sum_{i=1}^{n}g(R(u,e_i)X , R(u,e_i)Y)
                  =2\rho(X,Y)-\frac{1}{2}\, \alpha\,g(X,Y)
                   -\frac{1}{2}\,\beta\,g(X,u)g(Y,u).
\end{eqnarray}
\end{thm}

\medskip

{\noindent\bf Proof of Theorem \ref{thm4}.}  Let $T_1M =(T_1M, \eta,
\bar g, \phi, \xi)$ be the unit tangent bundle equipped with the
standard contact metric structure $(\eta, \bar g, \phi, \xi)$ and
assume that $T_1M$ is an $\eta$-Einstein manifold. Then, by the
definition, the Ricci tensor $\bar\rho$ of $T_1M$ takes of the
following form:
\begin{equation}\label{4-12}
\bar\rho=\alpha \bar g +\beta \eta\otimes\eta
\end{equation}
for some smooth functions $\alpha$ and $\beta$ on $T_1M$.

For a while, we adopt the traditional convention for the notations
in the classical tensor analysis. In the local coordinate
neighborhood, from (\ref{4-12}), we get
\begin{equation}\label{4-13}
\bar\rho_{ij}=\alpha \bar g_{ij} +\beta \eta_i\eta_j.
\end{equation}
Operating $\bar\nabla^i =\bar g^{ia}\bar\nabla_a$ on both sides of
(\ref{4-13}), we get
\begin{equation}\label{4-14}
\begin{aligned}
\bar\nabla^i \bar\rho_{ij}&=(\bar\nabla^i \alpha)\bar g_{ij} +(\bar\nabla^i \beta)\eta_i\eta_j
+\beta(\bar\nabla^i \eta_i)\eta_j +\beta\eta_i(\bar\nabla^i \eta_j)\\
&=\bar\nabla_j \alpha +(\bar\nabla_i \beta)\xi^i\eta_j +\beta(div\xi)\eta_j
+\beta\xi^i \bar\nabla_i \eta_j.
\end{aligned}
\end{equation}
Transvecting $\xi^j$ with (\ref{4-14}), we have
\begin{equation}\label{4-15}
\xi^j\bar\nabla^i \bar\rho_{ij}=\xi\alpha +\xi\beta +\beta(div\xi)+\beta(\bar\nabla_\xi \eta)\xi.
\end{equation}
Here, taking account of the second Bianchi identity, we get
$$
\bar\nabla^i \bar\rho_{ij} =\frac{1}{2} \bar\nabla_j \bar\tau
$$
and hence the left-hand side of (\ref{4-15}) implies $\frac{1}{2}\xi\bar\tau$.
Thus, from (\ref{4-10}), (\ref{4-11}) and (\ref{4-15})
we have
\begin{equation}\label{4-16}
\xi\bar\tau=2\xi\alpha +2\xi\beta.
\end{equation}
On one hand, by (\ref{4-13}), we get
$$
\bar\tau=(2n-1)\alpha +\beta.
$$
Thus, we have also
\begin{equation}\label{4-17}
\xi\bar\tau=(2n-1)\xi\alpha +\xi\beta.
\end{equation}
Then from (\ref{4-16}) and (\ref{4-17}), we have
\begin{equation}\label{4-18}
(2n-3)\xi\alpha -\xi\beta=0.
\end{equation}
Next, let $\bar X=(X^j)$ be a vector field on $T_1 M$ with
$\eta(\bar X)=0$. Transvecting $X^j$ with (\ref{4-14}), we have also
$$
X^j \bar\nabla^i \bar\rho_{ij}=\bar X\alpha +\beta(\bar \nabla_\xi \eta)(\bar X)
$$
and hence
\begin{equation}\label{4-19}
\frac{1}{2}\bar X\bar\tau =\bar X\alpha +\beta(\bar \nabla_\xi \eta)(\bar X).
\end{equation}
Here, we get
\begin{equation}\label{4-20}
\begin{aligned}
(\bar\nabla_\xi \eta)(\bar X)&=-\eta(\bar\nabla_\xi \bar X)\\
&=-\eta(\bar\nabla_{\bar X} \xi +[\xi,\bar X])\\
&=-\eta([\xi,\bar X]).
\end{aligned}
\end{equation}
On one hand, we get
\begin{equation}\label{4-21}
\begin{aligned}
-\eta([\xi,\bar X])&=\xi(\eta(\bar X))-\bar X(\eta(\xi))-\eta([\xi,\bar X])\\
&=d\eta(\xi, \bar X)\\
&=\bar g(\xi, \phi\bar X)\\
&=0.
\end{aligned}
\end{equation}
Thus from (\ref{4-19}) $\sim$ (\ref{4-21}), we have
\begin{equation}\label{4-22}
\bar X\bar\tau =2\bar X\alpha
\end{equation}
for vector field $\bar X$ with $\eta(\bar X)=0$. Since $\bar\tau=(2n-1)\alpha +\beta$
holds on $T_1M$, we have also
\begin{equation}\label{4-23-1}
\bar X\bar\tau =(2n-1)\bar X\alpha +\bar X\beta.
\end{equation}
Thus, by (\ref{4-22}) and (\ref{4-23-1}), we have
\begin{equation}\label{4-23}
(2n-3)\bar X\alpha +\bar X\beta=0
\end{equation}
for vector field $\bar X$ with $\eta(\bar X)=0$.

From now on, we state some fundamental properties of the
$\eta$-Einstein contact metric structure $(\eta, \bar g, \phi, \xi)$
on $T_1M$, by making use of the facts in the above. First of all, by
(\ref{4.2}), we see that the scalar curvature $\tau$ of the base
manifold $(M,g)$ (dim$M\geq2$) is constant.

Now setting $X=Y=e_j$ in (\ref{4.1}) and (\ref{4.3}) and taking sum for $j=1,\cdots,n$,
we obtain
\begin{eqnarray}
  \label{4-24} && \sum_{i,j=1}^{n}g(R(u,e_j)e_i , R(u,e_j)e_i)=(\alpha-4n+8)(n-1),\\
  \label{4-25} && \sum_{i,j=1}^{n}g(R(u,e_i)e_j , R(u,e_i)e_j)
                     =2\tau-\frac{1}{2}\,n\alpha -\frac{1}{2}\,\beta.
\end{eqnarray}
From (\ref{4-24}) and (\ref{4-25}), we have
\begin{equation}\label{4-26}
(3n-2)\alpha +\beta =4\tau +8(n-1)(n-2).
\end{equation}
Since $\tau$ is constant and $\xi=2u^h$, we have
\begin{equation}\label{4-27}
(3n-2)u^h \alpha +u^h\beta =0.
\end{equation}
And (\ref{4-18}) can be rewritten as follows :
\begin{equation}\label{4-28}
(2n-3)u^h \alpha -u^h\beta =0.
\end{equation}
From (\ref{4-27}) and (\ref{4-28}), we have
\begin{equation}\label{4-29}
u^h \alpha=0\;\; \text{and}\;\; u^h\beta =0.
\end{equation}
Operating $X^t\; (X\in T_pM) \in T_{(p,\;u)}T_1M$ on the both sides of (\ref{4-26}),
we have
\begin{equation}\label{4-30}
(3n-2)X^t \alpha +X^t\beta =0.
\end{equation}
Since $X^t$ is orthogonal to $\xi$ (i.e., $\eta(X^t)=0$), from (\ref{4-23}),
we have
\begin{equation}\label{4-31}
(2n-3)X^t \alpha +X^t\beta =0.
\end{equation}
Thus, from (\ref{4-30}) and (\ref{4-31}), we have
\begin{equation}\label{4-32}
X^t \alpha=0\;\; \text{and}\;\; X^t\beta =0\;\;\text{at}\;\;(p,u).
\end{equation}
Similarly, operating $X^h\; (X\in T_pM) \in T_{(p,\;u)}T_1M$
on the both sides of (\ref{4-26})
for vector field $X$ on $M$ such that $g(X,u)=0$, we have
\begin{equation}\label{4-33}
X^h \alpha=0\;\; \text{and}\;\; X^h\beta =0\;\;\text{at}\;\;(p,u).
\end{equation}
Summing up (\ref{4-29}), (\ref{4-32}) and (\ref{4-33}), we see that
the smooth functions $\alpha$ and $\beta$ are constants.
\hfill$\square$\medskip


\noindent By Theorem \ref{thm4}, we immediately obtain that
\begin{cor}\label{cor1}
$T_1M$ with $\eta$-Einstein structure has constant scalar curvature $\bar\tau$.
\end{cor}
\medskip


{\noindent\bf Proof of Theorem \ref{thm2}.} \quad
For $T_1M$ with constant scalar curvature it holds
\begin{equation}\label{4.6}
\sum_{i,j=1}^{n}g(R(u,e_j)e_i , R(u,e_j)e_i)=\frac{|R|^2}{n},
\end{equation}
where $|R|^2=\sum_{i,j,k=1}^{n}g(R(e_i,e_j)e_k , R(e_i,e_j)e_k)$
(\cite{BoeckxVan2}). From (\ref{4-24}), (\ref{4-25}) and
(\ref{4.6}), we have
\begin{eqnarray}
 \label{4.7} && \alpha =\frac{|R|^2}{n(n-1)} +4(n-2),\\
 \label{4.8} && \beta =4\tau -4n(n-2) -\frac{3n-2}{n(n-1)}|R|^2.
\end{eqnarray}
 Next, we integrate (\ref{4.3}) with $X=Y=u$ over $ S^{n-1}(1)$ in
$T_pM$. Then using the formula in \cite{ChenVanhecke}, we have
\begin{equation}\label{4.9}
  \frac1{n(n+2)}(|\rho|^2
  +\frac32|R|^2)=\frac{2\tau}{n}-\frac12\,\alpha-\frac12\,\beta.
\end{equation}
Eliminating $\alpha$ and $\beta$ from (\ref{4.7}) $\sim$ (\ref{4.9}), we
obtain the equation:
\begin{equation}\label{4.10}
  2|\rho|^2-3(n+1)|R|^2 =-4(n-1)(n+2)\tau +4n(n-1)(n-2)(n+2).
\end{equation}

\noindent In proof of Theorem \ref{thm4}, we know that $\alpha$,
$\beta$ and $\tau$ are constant. Since $\alpha$ is constant, from
(\ref{4.7}), we see that $|R|^2$ is constant. Thus, by (\ref{4.9}),
we see also that $|\rho|^2$ is constant. Therefore we have Theorem
\ref{thm2}.
\hfill$\square$\medskip


\begin{center}
\section{Special cases}\label{sec4}
\end{center}

(I) 2-dimensional case

It is well-known that $R(X,Y)Z=\kappa (g(Y,Z)X-g(X,Z)Y)$ always
holds. So, we have
$|R|^2=4\kappa^2,\;|\rho|^2=2\kappa^2,\;\tau=2\kappa$, where
$\kappa$ is the Gaussian curvature. From (\ref{4.10}), we see that
$M$ has Gaussian curvature $\kappa=0$ or $\kappa=1$. \vskip.3cm

(II) 3-dimensional case

 It is well-known that the curvature tensor $R$ of 3-dimensional
Riemannian manifold $(M,g)$ is of the following form.
\begin{equation}
\begin{aligned}\label{4-1}
  &R(X,Y,Z,W)=g(R(X,Y)Z,W)  \\
    &   =\Big\{g(X,W)\rho(Y,Z)+g(Y,Z)\rho(X,W)  \\
    &       \qquad -g(X,Z)\rho(Y,W)-g(Y,W)\rho(X,Z)\Big\}  \\
    &   \quad +\frac{\tau}{2}\{g(X,Z)g(Y,W)-g(Y,Z)g(X,W)\}
\end{aligned}
\end{equation}
for all vector fields $X$, $Y$, $Z$, $W$ on $M$. From (\ref{4-1}),
by direct calculation, we get
\begin{equation}\label{4-2}
|R|^2=4|\rho|^2 -\tau^2.
\end{equation}
By (\ref{4.10}) and (\ref{4-2}), we have
$$
23|\rho|^2 -6\tau^2 -20\tau +60=0,
$$
and thus
\begin{equation}\label{4-3}
23\left|\rho-\frac{\tau}{3} g \right|^2 +\frac{5}{3}(\tau-6)^2=0.
\end{equation}
From (\ref{4-3}), we have $\rho=\frac{\tau}{3} g$ and $\tau=6$ and hence
\begin{equation}\label{4-4}
\rho=2g.
\end{equation}
Thus by (\ref{4-1}) and (\ref{4-4}), we have
$$
R(X,Y,Z,W)=g(X,W)g(Y,Z)-g(X,Z)g(Y,W)
$$
and hence $(M,g)$ is a space of constant sectional curvature $1$.
The above result has been proved in \cite{BoeckxChoCuhn}. We may note
that our proof is much simpler than their proof. \vskip.3cm (III)
Conformally flat case

 By the similar arguments in \cite{BoeckxChoCuhn}, we can also
have the following.

\begin{thm}\label{thm3}
Let $M=(M,g)$ be an $n$-dimensional conformally flat manifold
$(n\geq4)$. Then $(T_1 M, \eta, \bar g, \phi, \xi)$ is
$\eta$-Einstein if and only if $(M,g)$ is a space of constant
sectional curvature $1$ or $n-2$.
\end{thm}
\vskip.3cm
\quad\;(IV) Einstein case

 Let $M=(M,g)$ be an
$n$-dimensional Einstein manifold $(n\geq3)$. Then we have
\begin{equation}\label{4-5}
\left|R+\frac{\tau}{2n(n-1)}g\cwedge\;g\right|^2=|R|^2 -\frac{2\tau^2}{n(n-1)}\;,
\end{equation}
where $(h\cwedge\;k)(X,Y,Z,W)=h(X,Z)k(Y,W)+h(Y,W)k(X,Z)-h(X,W)k(Y,Z)-h(Y,Z)k(X,W)$
for any (0,2)-tensors $h$ and $k$.
By (\ref{4.10}) and (\ref{4-5}), we have
\begin{equation*}
\begin{aligned}
\frac{2\tau^2}{n} &-3(n+1)\left\{\left|R+\frac{\tau}{2n(n-1)}g\cwedge\;g\right|^2
+\frac{2\tau^2}{n(n-1)}\right\}\\
&=-4(n-1)(n+2)\tau +4n(n-1)(n-2)(n+2)
\end{aligned}
\end{equation*}
and hence
\begin{equation}\label{4-6}
\begin{aligned}
-3&(n+1)\left|R+\frac{\tau}{2n(n-1)}g\cwedge\;g\right|^2\\
&=\frac{4(n+2)}{n(n-1)}\{\tau^2-n(n-1)^2\tau +n^2(n-1)^2(n-2)\}\\
&=\frac{4(n+2)}{n(n-1)}(\tau-n(n-1))(\tau-n(n-1)(n-2)).
\end{aligned}
\end{equation}
Then from (\ref{4-6}), we have
\begin{equation}\label{4-7}
n(n-1)\leq \tau \leq n(n-1)(n-2), \qquad n\geq3.
\end{equation}
By Theorem \ref{thm1}, we see that $(M,g)$ is super-Einstein by
virtue of (\ref{4.3}). Since the scalar curvature of $T_1M$ is
constant as shown by Theorem \ref{thm2}, this also follows from the result of
Boeckx and Vanhecke (\cite{BoeckxVan2}, Proposition 3.6.). Thus we
have

\begin{thm}\label{thm8}
Let $(M,g)$ be an $n$-dimensional Einstein manifold  and $(T_1 M,
\eta, \bar g, \phi, \xi)$ be the unit tangent bundle of $M$ equipped
with the canonical contact metric structure. If $T_1 M$ is
$\eta$-Einstein, then $M$ is super-Einstein and the scalar curvature
$\tau$ satisfies the above inequality (\ref{4-7}).
\end{thm}

In the remainder of this section, we shall consider the case that
$(M,g)$ is a $4$-dimensional Einstein manifold.

\medskip


{\noindent\bf Proof of Theorem \ref{thm7}.} \quad We may choose an
orthonormal basis $\{e_i\}$ (known as the Singer-Thorpe basis)
at each point $p\in M$ such that

\begin{equation}\label{5.1}
\begin{cases}
R_{1212}=R_{3434}=a,\quad R_{1313}=R_{2424}=b,\quad R_{1414}=R_{2323}=c,\\
R_{1234}=d,\quad R_{1342}=e,\quad R_{1423}=f,\\
R_{ijkl}=0 \;\; \text{whenever just three of the indices $i, j, k, l$ are distinct (\cite{ST})}.
\end{cases}
\end{equation}
Note that $d+e+f=0$ by the first Bianchi identity and
\begin{equation}\label{5.2}
a+b+c=-\frac{\tau}{4}.
\end{equation}
Further, by the direct calculation, we have
\begin{equation}\label{5.3}
\begin{aligned}
|R|^2 &=8(a^2 +b^2 +c^2 +d^2 +e^2 +f^2),\\
|\rho|^2 &=4(a+b+c)^2.
\end{aligned}
\end{equation}
From Theorem \ref{thm8}, since $M$ is super-Einstein, we have
(\cite{SV1, SV2})
\begin{equation}\label{5.4}
\pm d =a+\frac{\tau}{12} ,\quad \pm e =b+\frac{\tau}{12} ,\quad \pm f =c+\frac{\tau}{12} .
\end{equation}
From (\ref{4.1}), taking account of (\ref{5.1}), we have easily
\begin{eqnarray}
 \label{5.5} && 2(a^2 +d^2)=\alpha-8,\\
 \label{5.6} && 2(b^2 +e^2)=\alpha-8,\\
 \label{5.7} && 2(c^2 +f^2)=\alpha-8.
\end{eqnarray}
Thus from (\ref{5.5}) and (\ref{5.6}), taking account of
(\ref{5.4}), we have
\begin{equation}\label{5.8}
(a-b)(a+b+\frac{\tau}{12})=0.
\end{equation}
Similarly, we have
\begin{eqnarray}
 \label{5.9}  && (b-c)(b+c+\frac{\tau}{12})=0,\\
 \label{5.10} && (c-a)(c+a+\frac{\tau}{12})=0.
\end{eqnarray}
We first suppose that $a\neq b,\;b\neq c,\;c\neq a$. Then by
(\ref{5.8}) $\sim$ (\ref{5.10}), we get
\begin{equation*}
a+b+\frac{\tau}{12}=0,\quad b+c+\frac{\tau}{12}=0,\quad c+a+\frac{\tau}{12}=0.
\end{equation*}
Thus we have $a=b=c=-\frac{\tau}{24}$. However this is a
contradiction.

Next, we suppose that $a\neq b,\;b\neq c,\;c=a$ (i.e., $a=c,\; a\neq b$).
Then by (\ref{5.8}) $\sim$ (\ref{5.10}), we have
\begin{equation}\label{5.11}
a+b+\frac{\tau}{12}=0, \quad b+c+\frac{\tau}{12}=0.
\end{equation}
By (\ref{5.2}) and the hypothesis $a=c$, we have
\begin{equation}\label{5.12}
2a+b=-\frac{\tau}{4}.
\end{equation}
Thus by (\ref{5.11}) and (\ref{5.12}), we have
\begin{equation}\label{5.13}
a=c=-\frac{\tau}{6}, \quad b=\frac{\tau}{12}.
\end{equation}
Thus by (\ref{5.4}) and (\ref{5.13}), we have
\begin{equation}\label{5.14}
\pm d =-\frac{\tau}{12} ,\quad \pm e =\frac{\tau}{6} ,\quad \pm f =-\frac{\tau}{12} .
\end{equation}
Thus from (\ref{5.3}), (\ref{5.13}) and (\ref{5.14}), we have
\begin{equation}\label{5.15}
|R|^2 =\frac{5}{6} \tau^2,\quad |\rho|^2 =\frac{\tau^2}{4}.
\end{equation}
Then by (\ref{4.10}) and (\ref{5.15}), we obtain
\begin{equation}\label{5.16}
\tau^2 -6\tau +48=0.
\end{equation}
However, this quadratic equation (\ref{5.16}) does not admit a real
solution. This is also a contradiction. By the similar way, we can
also deduce a contradiction in the cases $b=c\neq a$ and $a=b\neq
c$. Thus, it must follow that $a=b=c$ and hence by (\ref{5.2}) and
(\ref{5.4}), we have
$$
a=b=c=-\frac{\tau}{12},\quad d=e=f=0.
$$
Therefore, by (\ref{5.1}), $(M,g)$ is a space of constant sectional curvature
$\frac{\tau}{12}$. Then we have
$$
|R|^2 =\frac{\tau^2}{6},\quad |\rho|^2 =\frac{\tau^2}{4}.
$$
Thus, by (\ref{4.10}), we have
\begin{equation}\label{5.17}
(\tau-12)(\tau-24)=0.
\end{equation}
Therefore, we have Theorem \ref{thm7}.
\hfill$\square$\medskip

\begin{center}{\bf Acknowledgments}
\end{center}Research of Y.D. Chai was partially supported by KOSEF
R01-2004-000-10183-0(2006), and research of J. H. Park was partially
supported by the Korea Research Foundation Grant funded by the
Korean Government (MOEHRD) KRF-2007-531-C00008.


\medskip

\noindent Y. D. Chai, S. H. ~Chun, J. H. ~Park :
 Department of Mathematics,
 Sungkyunkwan University,
 Suwon 440-746, Korea,
e-mail:  ydchai@skku.edu, cshyang@chonnam.ac.kr, parkj@skku.edu;
K.~Sekigawa :
    Department of Mathematics, Faculty of Science,
    Niigata University,
    Niigata, 950-2181, Japan,
    e-mail:  sekigawa@math.sc.niigata-u.ac.jp


\begin{thebibliography}{20}

\def\refprep#1#2{#1, {\em #2,} preprint.}
\def\refnot#1#2#3{#1, #2, {\em #3},\ to appear.}
\def\refart#1#2#3#4{#1, #2, {\em #3\/} #4.}
\def\refbook#1#2#3{#1, {\em #2,\/} #3.}


\bibitem{Besse}
Besse A (1978) Manifolds all of whose geodesics are closed.
Ergeb Math Grenzgeb 93, Berlin Heidelberg New-York: Springer

\bibitem{Blair1}
Blair DE (2002) Riemannian geometry of contact and symplectic
    manifolds. Progress in Math 203, Boston: Birkh\"auser


\bibitem{BoeckxVan}
Boeckx E, Vanhecke L (1997) Characteristic reflections on
    unit tangent sphere bundles. Houston J Math {\bf 23}: 427--448

\bibitem{BoeckxVan2}
Boeckx E, Vanhecke L (2001) Unit tangent sphere bundles
with constant scalar curvature. Czechoslovak Math J {\bf 51}: 523--544

\bibitem{BoeckxChoCuhn}
Boeckx E, Cho JT, Chun SH (2007) $\eta$-Einstein unit
tangent bundles. (preprint)

\bibitem{ChenVanhecke}
Chen B-Y, Vanhecke L (1981) Differential geometry of
geodesic spheres. J Reine Angew Math {\bf 325}: 28--67

\bibitem{Dombrowski}
Dombrowski P (1962) On the geometry of the tangent bundle.
J Reine Angew Math {\bf 210}: 73--88

\bibitem{SV1}
Sekigawa K, Vanhecke L (1986)
Volume-preserving geodesic symmetries on four-dimensional K\"ahler manifolds.
Lecture Notes in Math {\bf 1209} pp 275--291.
 Berlin: Springer

\bibitem{SV2}
Sekigawa K, Vanhecke L (1986) Volume-preserving geodesic
symmetries on four dimensional 2-stein spaces. Kodai math J
{\bf 9}: 215--224

\bibitem{ST}
Singer IM, Thorpe JA (1969) The curvature of $4$-dimensional
Einstein spaces. In: Global Analysis, Papers in Honor of K. Kodaira,
 pp 355--365: Princeton University Press

\bibitem{YanoKon}
Yano K, Kon M (1984) Structures on manifolds. Series in
Pure Mathematics Vol 3. Singapore: World Scientific


\end{thebibliography}
\end{document}